\documentclass{llncs}
\usepackage{amssymb, amsmath}
\usepackage{tikz}
\newcommand{\NP}{{\sf NP}}

\usepackage{comment}
\usepackage{seqsplit}
\title{New Bounds for the Snake-in-the-Box Problem}
\author{David Allison
\and
Dani\"el Paulusma}

\institute{
School of Engineering and  Computing Sciences, Durham University,\\
Durham, United Kingdom\\
\texttt{davidallison184@gmail.com}, \texttt{daniel.paulusma@durham.ac.uk}
}

\begin{document}
\maketitle
\setcounter{footnote}{0}

\begin{abstract}
The Snake-in-the-Box problem is that of finding a longest induced path in an $n$-dimensional hypercube.
We prove new lower bounds for the values $n\in \{11,12,13\}$.
The Coil-in-the-Box problem is that of finding a longest induced cycle in an $n$-dimensional hypercube.
We prove new lower bounds for the values $n\in \{10,11,12,13\}$.
\end{abstract}

\section{Introduction}\label{s-intro}

In 1958, Kautz~\cite{Ka58} observed that the vertices of an induced path in a hypercube can be used as code words in
Gray codes detecting single-bit errors. Applications of such codes date back to the 19th century, where they were used in telegraphy. They are nowadays used in many areas of mathematics, computer science and engineering, such as error correction in digital communication,  disk sector encoding, clock domain crossing, computer network topologies and in the design of genetic algorithms, just to name a few; we refer to Drapela~\cite{Dr15} for background information. The effectiveness of a Gray code depends on its length. Hence, the central research question is:

\medskip
\noindent
{\it How can we find a longest induced path in a hypercube?}
 
\medskip
\noindent
Kautz coined the name ``Snake-in-the-Box'' for this problem. The Snake-in-the-Box problem has been extensively studied in the literature, just as its variant ``Coil-in-the-Box'', where the goal is to determine a longest induced cycle instead of a path. As we shall discuss, the lengths of a longest snake and a longest coil are not known for any hypercube of dimension larger than~8.  The problem of finding a longest induced path is \NP-hard even for bipartite graphs~\cite{GJ79}, but its complexity for hypercubes is still unknown. As exhaustive search is too time consuming, a variety of computational and mathematical techniques have been used to find lower bounds for these problems. 
See Drapela~\cite{Dr15} for details on these techniques. 

In this note the main aim is to show three new lower bounds for Snake-in-the-Box and four new lower 
bounds for Coil-in-the-Box. We do this in Section~\ref{s-new} after first giving the required terminology in Section~\ref{s-terminology}. Examples meeting our lower bounds are presented in the Appendix. 
Our methods are based on the stochastic beam search method, which is a general heuristic search method that was recently used by Meyerson et al.~\cite{MDWP15,MWDP14} to improve lower bounds for
both Snake-in-the-Box and Coil-in-the-Box. We will explain our modifications to their algorithm in full detail in a future paper.

\section{Terminology}\label{s-terminology}
A graph $H$ is an {\it induced subgraph} of another graph $G$ if $H$ can be obtained from~$G$ by a sequence of vertex deletions. For $r\geq 1$, the {\it path} $P_r$ on $r$ vertices is the graph with vertices  $u_1,\ldots,u_r$ and edges $u_1u_2,\ldots,u_{r-1}u_r$. The {\it cycle} $C_r$ is the graph obtained from $P_r$ by adding the edge $u_ru_1$.
We say that $P_r$ and $C_r$ are of {\it length} $r-1$ and $r$, respectively. For $n\geq 1$, the $n$-dimensional hypercube $Q_n$ is the graph that contains $2^n$ vertices, each of which is represented by a binary vector of length $n$, such that  two vertices are adjacent if and only if their binary vectors differ by exactly one bit. 
A {\it snake} is an induced subgraph of $Q_n$ that is a path, and a {\it coil} is an induced subgraph of $Q_n$ that is a cycle.
The {\sc Snake-in-the-Box} problem is that of finding a longest snake in $Q_n$ for a given integer~$n$ and  
the {\sc Coil-in-the-Box} problem  is that of finding a longest coil in $Q_n$ for a given integer~$n$. 
Both problems have also been investigated for other graph classes (see, for example,~\cite{Wo98}).

\section{Known Results}\label{s-new}

Besides snakes and coils, a number of variants and generalizations, such as symmetric coils~\cite{AAC73} (also known as doubled coils), $k$-snakes and $k$-coils~\cite{HRSW15,Si66}, and single-track circuit codes~\cite{HPB96}  have been studied. However, in this section we will  only focus on the original notions of snakes and coils and we restrict our overview to exact or lower bounds for {\sc Snake-in-the-Box} and {\sc Coil-in-the-Box} for hypercubes of small dimensions, that is, we will consider $n$-dimensional hypercubes for $n\leq 20$. We refer to, for example,~\cite{EL00,Sn94,Ze97} for upper bounds for these two problems.

\medskip
\noindent
For $n\leq 8$, both {\sc Snake-in-the-Box} and {\sc Coil-in-the-Box} have been solved exactly:
\begin{itemize}
\item For $n\leq 5$, Kautz~\cite{Ka58} determined a longest coil. Using exhaustive search, Davies~\cite{Da65} found a longest snake for $n\leq 6$ and a longest coil for $n=6$. 
\item For $n=7$, Potter, Robinson, Miller and Kochut~\cite{PRMK94} used a genetic algorithm to find a longest snake and Kochut~\cite{Ko96} used exhaustive search to find a longest coil.
\item For $n=8$,  \"{O}sterg{\aa}rd and Pettersson used canonical augmentation to find a longest snake~\cite{OP15} and a longest coil~\cite{OP14}.
\end{itemize}
For $n\geq 9$, only lower bounds for the length of a longest snake or coil are known. These lower bounds have been improved over the years and the state of the art results are as follows:
\begin{itemize}
\item For $n=9$, Wynn~\cite{Wy12} proved that the length of a longest snake and coil is at least 190 and 188, respectively.
\item For $n=10$, Kinny~\cite{Ki12} showed that a longest snake has length at least~370 and Meyerson et al.~\cite{MDWP15} showed that a longest coil has length at least~362.
\item For $11\leq n\leq 12$, the results of Meyerson et al.~\cite{MDWP15} from 2015 are known to give the best lower bounds for the lengths of a longest snake and coil.
\item For $n=13$, Meyerson et al.~\cite{MDWP15} showed that a longest snake has length at least~2520 and Abbott and Katchalski~\cite{AK91} showed that a longest coil has length at least~2468.
\item For $14\leq n\leq 20$, the results of Abbott and Katchalski~\cite{AK91} from 1991 are still unbeaten. The same authors also showed a general lower bound by showing that for each $n\geq 21$, a longest coil has length at least $\frac{77}{256}2^n$. 
\end{itemize}

\noindent
We refer to Table~\ref{t-survey} for an overview of the currently best lower bounds on the length of a longest snake or coil for every $n\leq 20$. We note that the Snake-in-the-Box Records page http://ai1.ai.uga.edu/sib/sibwiki/ at the University of Georgia maintains a list of records for the Snake-in-the-Box and Coil-in-the-Box problems for every $n\leq 13$.

\begin{table}
\begin{center}
 \begin{small}
\begin{tabular}{|c|l|l|}
 \hline
Dimension & $\;\;$Snake Length & $\;\;$Coil Length\\
\hline
1 & $\;$ $1^*$ & $\;$ $0^*$ \\
\hline
2 &$\;$ 2$^*$ $\;$~\cite{Da65} &$\;$ 4$^*$ $\;$~\cite{Ka58}\\
\hline
3 &$\;$ 4$^*$ $\;$~\cite{Da65}  &$\;$  $6^*$ $\;$~\cite{Ka58}\\
\hline
4 & $\;$ 7$^*$ $\;$~\cite{Da65} & $\;$ $8^*$ $\;$~\cite{Ka58}\\
\hline
5 & $\;\;$13$^*$ $\;$\cite{Da65}  & $\;$ $14^*$~\cite{Ka58}\\
\hline
6 & $\;$ $26^*$$\;$~\cite{Da65} & $\;$ $26^*$~\cite{Da65}\\
\hline
7 & $\;$ 50$^*$~\cite{PRMK94} & $\;$ 48$^*$~\cite{Ko96}\\
\hline
8 & $\;$ 98$^*$~\cite{OP15} & $\;$ 96$^*$~\cite{OP14}\\
\hline
9 & $\;$ 190~\cite{Wy12} &$\;$ 188~\cite{Wy12}\\
\hline
10 &$\;$ 370~\cite{Ki12} & 
$\;$ {\bf 366} (362~\cite{MDWP15})\\
\hline
11 &$\;$ {\bf 712} $\;$ (707~\cite{MDWP15}) &$\;$ {\bf 692} (668~\cite{MDWP15})\\
\hline
12 & $\;$ {\bf 1373} (1302~\cite{MDWP15}) &$\;$ {\bf 1344} (1276~\cite{MDWP15})\\
\hline
13 &$\;$ {\bf 2687} (2520~\cite{MDWP15}) &$\;$ {\bf 2594} (2468~\cite{AK91})\\
\hline
14 &$\;$ 4932$\;\;\;$~\cite{AK91} &$\;$ 4934$\;\;$~\cite{AK91}\\
\hline
15 &$\;$ 9866$\;\;\;$~\cite{AK91} &$\;$ 9868$\;\;$~\cite{AK91}\\
\hline
16 &$\;$ 19738$\;$~\cite{AK91} &$\;$ 19740~\cite{AK91}\\
\hline
17 &$\;$ 39478$\;$~\cite{AK91} &$\;$ 39480~\cite{AK91}\\
\hline
18 &$\;$ 78958$\;$~\cite{AK91} &$\;$ 78960~\cite{AK91}\\
\hline
19 &$\;$ 157898~\cite{AK91} &$\;$ 157900~\cite{AK91}\\
\hline
20 &$\;$ 315798~\cite{AK91} &$\;$ 315800~\cite{AK91}\\
\hline
\end{tabular}
\end{small}
\end{center}
\caption{The lower bounds on the maximum length of a snake or coil in a hypercube of dimension $n=1,\ldots,20$. A * indicates that the bound is optimal. The unreferenced results in bold are the new bounds proven in this paper; the previous records are placed between parentheses after our lower bound values.
Note that the lower bounds on the provided snakes for $n\geq 14$ are deduced from the lower bound on the length of a longest coil.}\label{t-survey}
\end{table}

\section{Our Results}
We prove that the length of a longest snake is at least 712 for $n=11$, at least 1373 for $n=12$ and at least 2687 for $n=13$ and that the length of a longest coil is at least 
366 for $n=10$, at least 692 for $n=11$, at least 1344 for $n=12$ and at least 2594 for $n=13$; see also Table~\ref{t-survey}. We do this by giving examples of snakes and coils of these lengths in the Appendix.  We checked the correctness of our solutions using the verifier provided at \url{http://ai1.ai.uga.edu/sib/sibwiki/doku.php/checker}. 

\section*{Acknowledgements}
This work made use of the facilities of the Hamilton HPC Service of Durham University.

\bibliographystyle{plain}

\appendix

\section{The New Snakes and Coils}

Recall that two vertices in a hypercube are adjacent if and only if they differ in exactly one bit. 
A path $u_1u_2\cdots u_n$ in a hypercube $Q_n$ can be represented by a {\it transition sequence} $i_1,\ldots,i_{n-1}$ where for $j=1,\ldots,n-1$, $i_j$  denotes the vector entry of the bit in which vertices $u_{j-1}$ and $u_j$ differ. Similarly, a cycle in a hypercube can be represented by a transition sequence as well.

\subsection{Snake of length 712 (Dimension: 11)}
\seqsplit{0,1,2,3,0,1,4,0,3,5,4,0,1,4,5,2,3,5,4,1,0,4,6,3,5,0,3,4,1,3,2,1,0,3,5,0,1,4,5,0,3,5,4,2,3,7,1,5,3,1,2,3,0,1,4,0,3,5,4,0,1,4,5,2,3,5,1,6,5,0,3,5,4,1,2,4,5,3,0,5,4,1,0,5,3,0,1,2,3,0,8,2,0,1,4,0,3,5,4,0,1,4,5,2,3,5,4,1,0,4,6,5,0,4,1,3,2,1,0,3,5,0,1,4,5,0,3,5,4,2,1,7,2,4,5,3,0,5,2,4,1,3,2,4,5,2,0,4,2,3,5,4,2,6,4,3,2,4,5,3,0,5,4,1,0,5,3,0,1,2,3,1,4,0,3,5,0,8,9,0,5,3,0,4,1,3,2,1,0,3,5,0,1,4,5,0,3,5,4,2,3,4,6,2,4,5,3,2,4,0,2,5,4,2,3,1,4,2,5,0,3,5,4,2,7,1,2,4,5,3,0,5,4,1,0,5,3,0,1,2,3,1,4,0,5,6,4,0,1,3,2,1,0,4,5,3,0,4,1,0,3,2,1,0,4,1,3,0,8,2,0,5,4,1,0,5,3,0,1,2,3,1,4,5,6,1,5,3,2,5,4,1,0,4,5,3,0,4,1,0,3,2,1,3,5,1,7,3,2,4,5,3,0,5,4,1,0,5,3,0,1,2,3,1,4,3,0,5,3,6,4,0,1,4,5,3,2,5,4,1,0,4,5,3,0,4,1,0,3,2,1,0,10,8,0,1,2,3,0,1,4,0,3,5,4,0,1,4,5,2,3,5,4,1,0,4,6,3,5,0,3,4,1,3,2,1,0,3,5,0,1,4,5,0,3,5,4,2,3,7,1,5,3,1,2,3,0,1,4,0,3,5,4,0,1,4,5,2,3,5,1,6,5,0,3,5,4,1,2,4,5,3,0,5,4,1,0,5,3,0,1,2,3,0,8,2,0,1,4,0,3,5,4,0,1,4,5,2,3,5,4,1,0,4,6,5,0,4,1,3,2,1,0,3,5,0,1,4,5,0,3,5,4,2,1,7,2,4,5,3,0,5,2,4,1,3,2,4,5,2,0,4,2,3,5,4,2,6,4,3,2,4,5,3,0,5,4,1,0,5,2,4,3,5,0,7,6,9,3,1,5,3,2,5,4,1,0,4,5,3,0,4,1,0,3,2,1,3,5,1,7,3,2,4,5,3,0,5,4,1,0,5,3,0,1,2,3,1,4,3,0,5,3,6,4,0,1,4,5,3,2,5,4,1,0,4,5,3,0,4,1,0,2,8,0,3,2,1,0,3,5,0,1,4,5,0,3,5,4,2,1,4,5,0,3,5,6,4,0,5,3,1,2,3,0,1,4,3,7,4,1,2,3,1,4,0,3,5,4,1,2,4,5,3,0,5,4,1,0,4,7,3,4,0,1,4,2,5,3,0,5,6,7,0,4,1,0,5,3,1,4,0,2,5,4,1,0,4,5,3,0,4,1,0,5,8,1,5,3,0,4,5,3,1,5,0,2,1,0,4,1,3,5,1,0,4,1,10,9,1,4,0}
\subsection{Snake of length 1373 (Dimension: 12)}
\seqsplit{0,1,2,3,4,5,2,3,1,0,3,2,5,3,1,2,3,4,5,2,3,1,0,3,2,6,1,2,3,4,1,3,5,4,3,0,1,3,4,5,2,3,4,1,7,6,1,4,3,2,5,4,3,1,0,3,4,5,3,1,4,3,2,5,4,3,1,0,6,1,2,5,4,3,2,5,1,2,4,5,2,0,1,2,5,4,3,2,5,8,0,4,3,2,5,4,3,1,2,6,0,1,2,3,4,5,2,3,1,0,3,2,5,3,1,2,3,4,5,2,3,1,0,3,2,6,0,7,1,2,3,1,0,2,6,0,1,4,3,2,5,4,3,1,0,3,4,5,3,1,4,3,2,5,4,3,1,0,6,1,4,0,2,3,4,5,2,1,5,4,9,5,6,2,3,0,1,3,2,5,4,3,2,1,3,5,2,3,0,1,3,2,5,4,3,2,1,0,6,2,1,3,4,5,2,3,4,0,8,5,2,3,4,5,2,1,0,2,5,4,2,1,5,2,3,4,5,2,1,6,0,1,3,4,5,2,3,4,1,3,5,4,3,0,1,3,4,5,2,3,4,1,0,7,8,0,1,4,3,2,5,4,3,1,0,3,4,5,3,1,4,3,2,5,4,3,1,0,6,1,2,5,4,3,2,5,1,2,4,5,2,0,1,2,5,4,3,2,1,8,2,3,4,2,1,3,2,5,4,3,2,1,0,2,3,6,1,4,5,1,3,5,2,0,1,2,5,4,3,2,0,6,2,3,4,5,2,3,1,0,10,3,1,2,5,4,3,2,5,8,0,4,3,2,5,4,3,1,2,6,0,1,2,3,4,5,2,3,1,0,3,2,5,3,1,2,3,4,5,2,3,1,0,3,2,6,5,8,0,1,5,2,3,4,5,2,1,6,0,1,3,4,5,2,3,4,1,3,5,4,3,0,1,3,4,5,2,3,4,1,6,7,1,4,3,2,5,4,3,1,0,3,4,5,3,1,5,2,3,4,5,1,6,2,3,0,1,4,5,1,3,5,2,3,0,1,3,2,5,4,3,2,1,0,3,8,4,3,0,1,3,4,5,2,3,4,1,3,5,4,3,0,1,3,4,5,2,3,4,1,0,6,4,1,3,2,5,4,3,2,0,4,2,5,4,3,2,1,4,2,5,4,9,1,4,5,2,3,4,1,5,6,0,1,2,3,4,5,2,3,1,0,3,2,5,3,1,2,3,4,5,2,3,1,0,3,2,6,1,2,3,4,5,2,1,8,6,1,3,4,5,2,3,4,1,3,5,4,3,0,1,3,4,5,2,3,4,1,0,6,1,5,2,3,4,5,2,1,0,2,5,3,1,2,3,4,5,2,3,1,7,6,1,3,2,5,4,3,2,1,3,5,2,0,1,2,5,4,3,2,1,0,6,1,2,3,4,5,2,3,0,2,4,3,2,5,4,1,8,5,1,2,5,4,3,2,1,5,6,0,1,4,3,2,5,4,3,1,0,3,4,5,3,1,4,3,2,5,4,3,1,0,3,4,6,1,4,3,2,5,4,3,1,0,11,6,3,0,1,3,2,5,4,3,2,1,3,5,2,3,0,1,3,2,5,4,3,2,1,0,6,2,1,3,4,5,2,3,4,0,8,5,2,3,4,5,2,1,3,5,6,1,3,4,5,2,3,4,1,3,5,4,3,0,1,3,4,5,2,3,4,1,6,7,1,4,3,2,5,4,3,1,0,3,4,5,3,1,4,3,2,1,6,2,3,0,1,4,5,1,3,5,2,3,0,1,3,2,5,4,3,2,1,0,3,8,4,3,0,1,3,4,5,2,3,4,1,3,5,4,3,0,1,3,4,5,2,3,4,1,0,6,2,0,1,3,2,5,3,0,2,3,4,5,2,1,5,4,9,5,6,2,3,0,1,3,2,5,4,3,2,1,3,5,2,3,0,1,3,2,5,4,3,2,1,0,6,2,1,3,4,5,2,3,4,0,8,5,2,3,4,5,2,1,0,2,5,4,2,1,5,2,3,4,5,2,1,6,0,1,3,4,5,2,3,4,1,3,5,4,3,0,1,3,4,5,2,3,4,1,0,7,1,5,2,3,4,5,2,1,0,3,6,1,3,2,5,4,3,2,8,6,5,4,3,2,5,4,1,10,4,5,2,3,0,1,3,2,5,4,3,2,1,0,6,2,1,3,4,5,2,3,4,0,8,5,2,3,4,5,2,1,0,2,5,4,2,1,5,2,3,4,5,2,1,6,0,1,3,4,5,2,3,4,1,3,5,4,3,0,1,3,4,5,2,3,4,1,6,7,1,4,3,2,5,4,3,1,0,2,6,1,2,3,4,5,2,3,1,0,3,2,5,3,1,5,4,3,2,5,1,8,4,3,0,1,3,4,5,2,3,4,1,3,5,4,3,0,1,3,4,5,2,3,4,1,0,6,4,1,3,2,5,4,3,2,0,4,2,5,4,3,2,1,4,2,5,4,9,1,4,5,2,3,4,1,5,6,0,1,2,3,4,5,2,3,1,0,3,2,5,3,1,2,3,4,5,2,3,1,0,3,2,6,1,2,3,4,5,2,1,8,6,1,2,5,4,3,2,5,1,2,3,0,1,3,4,5,2,3,4,1,0,6,7,8,6,0,1,4,3,2,5,4,3,1,0,3,4,5,3,1,4,3,2,5,4,3,1,0,3,4,6,1,3,0,2,3,4,5,2,1,0,3,8,1,3,2,5,4,3,1,2,6,0,1,2,3,4,5,2,1,0,2,5,3,1,2,3,4,5,2,3,1,0,3,6,1,3,4,5,2,3,1,4,10,0,1,4,3,2,5,1,8,0,1,5,2,3,4,9,6,1,8,0,1,4,3,2,5,4,3,1,9,6,7,10,1,2,5,3,0,2,3,4,5,2,1,0,5,2,3,5,9,0,1,5,2,3,4,5,1,2,10,3,7,4,1,5,2,3,1,2,4,8,3,6,4,1,5,4,3,6,8,2,1,0,8,3,2,1,4,5,1,3,2,10,0,3,10,11,1,2,5,1,3,11,7,5,3,11,10,7,9,11,0,1,5,4,3,1,9,0,6,3,0,1,4,5,2,3,4,5}

\subsection{Snake of Length 2687 (Dimension: 13)}
\seqsplit{0,1,2,3,4,5,2,3,1,0,3,2,5,3,1,2,3,4,5,2,3,1,0,3,2,6,1,2,3,4,1,3,5,4,3,0,1,3,4,5,2,3,4,1,7,6,1,4,3,2,5,4,3,1,0,3,4,5,3,1,4,3,2,5,4,3,1,0,6,1,2,5,4,3,2,5,1,2,4,5,2,0,1,2,5,4,3,2,5,8,0,4,3,2,5,4,3,1,2,6,0,1,2,3,4,5,2,3,1,0,3,2,5,3,1,2,3,4,5,2,3,1,0,3,2,6,0,7,1,2,3,1,0,2,6,0,1,4,3,2,5,4,3,1,0,3,4,5,3,1,4,3,2,5,4,3,1,0,6,1,4,0,2,3,4,5,2,1,5,4,9,5,6,2,3,0,1,3,2,5,4,3,2,1,3,5,2,3,0,1,3,2,5,4,3,2,1,0,6,2,1,3,4,5,2,3,4,0,8,5,2,3,4,5,2,1,0,2,5,4,2,1,5,2,3,4,5,2,1,6,0,1,3,4,5,2,3,4,1,3,5,4,3,0,1,3,4,5,2,3,4,1,0,7,1,5,2,3,4,5,2,1,0,2,5,3,1,2,3,4,5,2,3,1,6,3,2,0,1,5,4,3,2,5,1,2,0,1,3,2,5,4,3,2,1,0,8,1,2,3,4,5,2,1,0,2,5,4,2,1,5,2,3,4,5,2,1,6,0,1,3,4,5,2,3,4,1,3,5,4,3,0,1,3,4,5,2,3,4,1,3,10,4,3,0,1,3,2,5,4,3,2,1,0,6,2,1,3,4,5,2,3,4,0,8,5,2,3,4,5,2,1,0,2,5,4,2,1,5,2,3,4,5,2,1,6,0,1,3,4,5,2,3,4,1,3,5,4,3,0,1,3,4,5,2,3,4,1,6,7,1,4,3,2,5,4,3,1,0,3,4,5,3,1,4,3,2,1,6,2,3,0,1,4,5,1,3,5,2,3,0,1,3,2,5,4,3,2,1,0,3,8,4,3,0,1,3,4,5,2,3,4,1,3,5,4,3,0,1,3,4,5,2,3,4,1,0,6,4,1,3,2,5,4,3,2,0,4,2,5,4,3,2,1,4,2,5,4,9,1,4,5,2,3,4,1,5,6,0,1,2,3,4,5,2,3,1,0,3,2,5,3,1,2,3,4,5,2,3,1,0,3,2,6,1,5,4,3,2,5,1,0,8,1,4,3,2,5,4,3,1,6,0,1,3,4,5,2,3,4,1,3,5,4,3,0,1,3,4,5,2,3,4,1,0,6,1,5,2,3,4,5,2,1,5,4,7,1,5,2,3,4,5,2,1,5,0,1,3,2,5,4,3,2,1,6,0,1,2,3,4,5,2,1,4,2,3,0,1,3,2,5,4,3,2,1,6,8,1,2,3,4,1,6,4,3,0,1,3,4,5,2,3,4,1,3,5,4,3,0,1,3,4,5,2,3,4,1,0,6,5,1,2,3,4,5,2,1,5,9,1,3,2,1,6,0,1,5,4,1,0,6,10,11,6,0,1,3,2,5,4,3,2,1,3,5,2,3,0,1,3,2,5,4,3,2,1,0,6,2,1,3,4,5,2,3,4,0,8,5,2,3,4,5,2,1,0,2,5,4,2,1,5,2,3,4,5,2,1,6,0,1,3,4,5,2,3,4,1,3,5,4,3,0,1,3,4,5,2,3,4,1,6,7,1,4,3,2,5,4,3,1,0,3,4,5,3,1,4,3,2,1,6,2,3,0,1,4,5,1,3,5,2,3,0,1,3,2,5,4,3,2,1,0,3,8,4,3,0,1,3,4,5,2,3,4,1,3,5,4,3,0,1,3,4,5,2,3,4,1,0,6,4,1,3,2,5,4,3,2,0,4,2,5,4,3,2,1,4,2,5,4,9,1,4,5,2,3,4,1,5,6,0,1,2,3,4,5,2,3,1,0,3,2,5,3,1,2,3,4,5,2,3,1,0,3,2,6,1,2,3,4,1,0,8,1,4,3,2,5,4,3,1,6,0,1,3,4,5,2,3,4,1,3,5,4,3,0,1,3,4,5,2,3,4,1,0,6,1,5,2,3,4,5,2,1,5,4,7,1,5,2,3,4,5,2,1,5,0,1,3,2,5,4,3,2,1,6,0,1,2,3,4,5,2,1,4,2,3,0,1,3,2,5,4,3,2,1,6,8,3,5,2,0,1,2,5,4,3,2,1,5,6,0,1,4,3,2,5,4,3,1,0,3,4,5,3,1,4,3,2,5,4,3,1,0,3,4,6,1,5,2,3,1,10,6,1,3,2,5,4,3,2,1,3,5,2,3,0,1,3,2,5,4,3,2,1,0,6,1,5,4,3,2,5,4,1,8,2,1,5,2,3,4,5,2,1,6,0,1,3,4,5,2,3,4,1,3,5,4,3,0,1,3,4,5,2,3,4,1,6,7,1,4,3,2,5,4,3,1,0,3,4,5,3,1,4,3,2,1,6,2,3,0,1,4,5,1,3,5,2,3,0,1,3,2,5,4,3,2,1,0,3,8,5,3,1,4,3,2,5,4,3,1,0,3,4,5,3,1,4,3,2,5,4,3,1,0,6,1,2,5,4,3,1,4,2,0,1,2,5,4,3,2,1,5,2,9,6,0,1,3,2,5,4,3,2,1,3,5,2,3,0,1,3,2,5,4,3,2,1,0,6,2,1,3,4,5,2,3,4,0,8,5,2,3,4,5,2,1,0,2,5,4,2,1,5,2,3,4,5,2,1,6,0,1,3,4,5,2,3,4,1,3,5,4,3,0,1,3,4,5,2,3,4,1,0,7,1,5,2,3,4,5,2,1,0,2,5,3,1,2,3,4,5,2,3,1,6,3,2,0,1,5,4,3,2,5,1,2,0,1,3,2,5,4,3,2,1,0,8,5,0,6,1,4,3,2,5,4,3,1,0,3,4,5,3,1,5,2,1,0,3,4,6,1,4,3,2,1,3,5,2,0,1,2,5,4,3,5,9,3,2,5,4,3,2,1,7,9,2,1,3,4,1,5,6,7,10,6,5,1,4,3,1,2,3,5,12,6,5,4,3,0,1,3,4,5,2,3,4,1,0,6,4,11,6,1,3,2,5,4,3,2,1,0,6,2,1,3,4,5,2,3,4,0,8,5,2,3,4,5,2,1,0,2,5,4,2,1,5,2,3,4,5,2,1,6,0,1,3,4,5,2,3,4,1,3,5,4,3,0,1,3,4,5,2,3,4,1,6,7,1,4,3,2,5,4,3,1,0,3,4,5,3,1,4,3,2,1,6,2,3,0,1,4,5,1,3,5,2,3,0,1,3,2,5,4,3,2,1,0,3,8,4,3,0,1,3,4,5,2,3,4,1,3,5,4,3,0,1,3,4,5,2,3,4,1,0,6,5,1,2,3,4,5,2,1,0,2,5,4,9,1,4,5,2,3,4,1,5,6,0,1,2,3,4,5,2,3,1,0,3,2,5,3,1,2,3,4,5,2,3,1,0,3,2,6,1,2,3,4,1,0,8,1,4,3,2,5,4,3,1,6,0,1,3,4,5,2,3,4,1,3,5,4,3,0,1,3,4,5,2,3,4,1,0,6,1,5,2,3,4,5,2,1,5,4,7,1,5,2,3,4,5,2,1,5,0,1,3,2,5,4,3,2,1,6,0,1,2,3,4,5,2,1,4,2,3,0,1,3,2,5,4,3,2,1,6,8,3,5,2,0,1,2,5,4,3,2,1,5,6,0,1,4,3,2,5,4,3,1,0,3,4,5,3,1,4,3,2,5,4,3,1,0,3,4,6,1,5,2,3,1,10,6,1,3,2,5,4,3,2,1,3,5,2,3,0,1,3,2,5,4,3,2,1,0,6,1,5,4,3,2,5,4,1,8,2,1,5,2,3,4,5,2,1,6,0,1,3,4,5,2,3,4,1,3,5,4,3,0,1,3,4,5,2,3,4,1,6,7,1,4,3,2,5,4,3,1,0,3,4,5,3,1,4,3,2,1,6,2,3,0,1,4,5,1,3,5,2,3,0,1,3,2,5,4,3,2,1,0,3,8,5,3,1,4,3,2,5,4,3,1,0,3,4,5,3,1,4,3,2,5,4,3,1,0,6,1,2,5,4,3,1,4,2,0,1,2,5,4,3,2,1,5,2,9,6,0,1,3,2,5,4,3,2,1,3,5,2,3,0,1,3,2,5,4,3,2,1,0,6,2,1,3,4,5,2,3,4,0,8,5,2,3,4,5,2,1,0,2,5,4,2,1,5,2,3,4,5,2,1,6,0,1,3,4,5,2,3,4,1,3,5,4,3,0,1,3,4,5,2,3,4,1,0,7,1,5,2,3,4,5,2,1,0,2,5,3,1,2,3,4,5,2,3,1,6,3,2,0,1,5,4,3,2,5,1,2,0,1,3,2,5,4,3,2,1,0,8,5,0,6,1,4,3,2,5,4,3,1,0,3,4,5,3,1,4,6,5,3,4,5,2,1,0,2,5,3,1,2,3,4,1,6,4,3,0,1,3,4,5,1,7,6,5,2,11,6,3,0,1,3,4,5,2,3,4,1,0,6,5,8,2,0,1,3,4,5,2,3,4,1,7,6,1,4,3,2,5,4,3,1,0,3,4,5,3,1,4,3,2,5,4,3,1,0,6,1,2,5,4,3,2,5,1,2,4,5,2,0,1,2,5,4,3,2,5,8,0,4,3,2,5,4,3,1,2,6,0,1,2,3,4,5,2,3,1,0,3,2,5,3,1,2,3,4,5,2,3,1,0,3,6,2,5,1,4,3,1,2,5,9,6,3,5,4,3,0,1,3,4,5,2,3,4,1,6,0,5,1,2,3,4,5,2,1,0,2,5,3,1,5,4,3,2,5,1,6,4,3,0,1,3,4,5,1,7,2,3,0,1,4,5,1,3,5,2,3,0,1,3,2,5,4,3,2,1,0,6,2,3,1,4,5,2,3,4,0,2,3,9,6,2,5,4,3,1,0,3,4,8,3,1,0,2,3,4,5,2,3,1,0,3,2,5,3,1,5,4,1,0,3,2,6,1,5,10,6,5,2,3,4,1,3,5,4,3,0,1,3,4,5,2,3,4,1,0,6,1,8,0,4,3,2,5,4,2,11,4,3,2,5,4,1,11,5,2,6,1,0,2,3,4,5,2,3,1,0,3,4,1,2,3,4,5,1,3,4,9,1,4,5,2,3,4,1,0,6,2,7,6,0,1,3,2,5,4,3,2,1,3,5,2,3,0,1,3,2,5,4,3,2,1,0,6,2,1,3,4,5,2,3,4,0,8,5,2,1,6,0,1,3,4,5,2,3,4,1,0,10,4,1,3,2,5,4,3,2,0,1,6,2,1,3,4,5,2,3,1,2,10,6,0,1,4,3,2,5,4,3,1,0,6,1,3,4,1,6,7,3,1,0,2,3,4,5,2,3,1,0,3,2,5,3,1,5,4,1,6,2,1,4,3,2,5,4,1,9,2,5,1,2,3,4,1,2,7,1,3,2,5,4,3,1,2,6,0,1,2,3,4,5,2,1,0,2,5,3,1,2,3,4,5,2,3,1,0,6,1,4,8,6,5,2,3,4,5,2,6,5,3,2,0,4,2,1,5,6,11,1,5,2,3,4,5,2,1,0,6,7,11,2,5,6,9,1,5,2,3,4,5,1,3,5,2,6,1,2,5,4,3,2,0,4,2,5,4,3,2,1,10,9,2,8,9,1,3,4,5,2,4,0,2,3,4,5,2,3,1,4,2,7,1,4,5,1,6,0,1,5,4,3,2,1,7,3,1,4,3,2,5,4,3,1,0,3,4,5,3,1,4,3,2,5,4,3,1,0,3,6,4,0,8,4,1,3,2,9,0,5,1,6,0,3,6,8,5,2,10,5,1,2,5,3,2,7,1,2,3,9,8,4,2,1,4,3,1,7,8,2,5,1,4,3,2,1,7,2,0,9,8,5,2,4,0,11,1,3,2,1,7,10,9,7,2,1,3,4,5,2,1,12,0,1,3}

\subsection{Coil of length 366 (Dimension: 10)}
\seqsplit{0,1,2,3,0,4,5,0,3,2,0,4,3,0,1,2,3,0,4,6,0,3,2,1,0,5,4,0,1,2,0,5,2,3,0,1,2,0,5,4,0,2,1,0,7,5,0,1,2,3,0,5,4,0,3,2,0,5,3,0,1,2,3,0,5,6,0,3,2,1,0,4,5,0,1,2,0,4,2,3,0,1,2,0,4,5,0,2,8,7,2,0,5,4,0,2,1,0,3,2,4,0,2,1,0,5,4,0,1,2,3,0,6,5,0,3,2,1,0,3,5,0,2,3,0,4,5,0,3,2,1,0,3,5,0,2,3,7,1,0,5,1,3,0,2,1,0,5,4,0,1,2,0,6,4,2,5,3,0,2,4,3,2,1,0,2,6,0,3,2,1,0,5,6,1,2,0,3,4,9,3,0,5,4,0,2,1,0,3,2,4,0,2,1,0,5,4,0,1,2,3,0,6,5,0,3,2,1,0,3,5,0,2,3,0,4,5,0,3,2,1,0,5,7,0,1,2,0,4,5,0,2,1,0,3,2,5,0,2,1,0,4,5,0,1,2,3,0,6,4,0,3,2,1,0,3,4,0,2,3,0,5,4,0,3,2,1,0,3,4,0,2,3,8,7,3,2,0,3,4,0,1,2,3,0,4,5,0,3,2,0,4,3,0,1,2,3,0,4,6,0,3,2,1,0,5,4,0,1,2,0,5,2,3,0,1,2,0,5,4,0,3,7,4,3,2,0,4,6,2,4,5,0,2,1,0,3,4,6,3,0,1,2,0,3,4,2,0,1,2,3,0,5,4,3,6,4,0,2,1,4,9,1,0,4}
\subsection{Coil of length 692 (Dimension: 11)}
\seqsplit{0,1,2,0,3,2,4,0,5,4,2,1,5,4,0,5,1,3,0,1,2,6,1,4,2,1,5,0,3,5,1,2,4,1,5,0,4,1,2,4,0,3,2,7,3,4,2,3,0,4,5,0,3,1,5,0,4,5,1,2,4,1,6,2,1,5,4,2,1,0,2,5,1,2,4,3,0,4,2,1,5,4,2,0,8,4,0,1,2,0,3,2,4,0,5,4,2,1,5,4,0,5,1,3,0,1,6,4,1,2,4,5,1,3,4,1,2,4,0,3,1,4,3,5,1,2,4,1,3,7,5,4,0,1,2,4,5,0,3,1,5,0,4,5,1,2,4,1,6,0,1,3,0,5,1,2,0,1,3,0,4,2,1,0,2,5,1,3,0,1,2,0,9,5,0,1,2,0,3,2,4,0,5,4,2,1,5,4,0,5,1,3,0,1,2,6,1,4,2,1,5,0,3,5,1,2,4,1,5,0,4,1,2,4,0,3,2,7,3,4,2,3,0,4,5,0,3,1,5,0,4,5,1,2,4,1,6,0,1,3,0,5,1,2,0,1,3,0,4,2,1,0,2,5,1,3,0,1,2,8,5,1,2,0,3,2,4,0,5,4,2,1,5,4,0,5,1,3,0,1,6,4,1,2,4,5,1,3,4,1,2,4,0,3,1,4,3,5,1,2,4,1,3,6,7,3,5,0,4,1,2,4,0,3,2,0,5,4,0,3,2,1,5,3,0,4,6,3,0,5,3,4,1,5,0,2,1,5,3,0,5,1,4,0,5,2,10,5,4,2,1,4,0,5,1,4,2,1,5,3,4,1,2,0,5,4,2,1,4,6,1,0,3,1,5,0,4,5,1,2,4,5,0,4,2,3,0,2,1,5,8,2,1,0,3,1,5,2,0,1,2,4,0,3,1,0,2,1,5,0,3,1,0,6,1,4,2,1,5,4,0,5,1,3,0,5,4,2,1,0,4,5,7,3,1,4,2,1,5,3,4,1,3,0,4,2,1,4,3,1,5,4,2,1,4,6,1,0,3,1,5,0,4,5,1,2,4,5,0,4,2,3,0,2,1,0,4,8,1,2,3,0,4,2,1,0,2,5,1,2,4,1,0,3,1,4,2,1,5,4,6,3,4,5,1,2,4,5,0,4,2,3,0,2,1,5,0,2,4,1,2,9,0,4,2,1,4,6,1,0,3,1,5,0,4,5,1,2,4,5,0,4,2,3,0,2,1,5,8,2,1,0,3,1,5,2,0,1,2,4,0,3,1,0,2,1,5,0,3,1,0,6,1,4,2,1,5,4,0,5,1,3,0,5,4,2,1,0,8,2,4,1,3,0,4,2,1,4,0,5,4,10,3,7,10,4,5,0,4,1,2,4,5,0,3,1,5,0,4,5,1,2,4,1,6,5,4,1,2,4,5,8,3,5,4,2,5,8,4,0,5,2,1,0,4,5,0,1,3,2,1,8,0,2,5,4,0,5,3,2,0,6,2,4,1,0,2,3,4,1,5,4,0,1,4,6,0,4,2,6,7,2,4,3,7,2,1,4,6,2,8}
\subsection{A Coil of length 1344 (Dimension: 12)}
\seqsplit{0,1,2,3,4,5,0,3,4,2,3,0,1,2,3,4,2,6,4,0,2,4,3,0,5,3,1,5,4,3,0,5,3,2,4,5,3,0,5,4,1,0,4,7,0,5,4,1,0,3,2,4,3,0,5,4,3,2,1,0,2,3,5,6,1,0,3,2,1,0,4,1,3,0,5,3,1,2,4,3,2,0,1,3,2,4,3,1,8,4,3,2,0,1,2,3,4,5,0,3,4,2,3,0,1,2,3,4,2,6,4,1,2,3,0,5,3,2,4,5,3,0,5,4,1,2,4,5,0,3,5,7,3,4,5,0,2,1,0,3,2,4,3,0,5,4,3,1,6,4,1,0,3,4,5,0,3,2,1,0,2,4,3,0,2,1,5,4,1,0,4,3,9,4,0,1,2,3,4,5,0,3,4,2,3,0,1,2,3,4,2,6,4,0,2,4,3,0,5,3,1,5,4,3,0,5,3,2,4,5,3,0,5,4,1,0,4,7,5,0,4,1,0,3,2,4,3,0,5,4,3,2,1,4,6,1,3,5,0,3,1,4,0,1,2,3,0,1,4,0,5,4,2,0,1,3,2,4,1,8,4,2,0,1,2,3,4,5,0,3,4,2,3,0,1,2,3,4,2,6,4,0,3,5,0,1,4,5,0,3,5,4,2,3,5,0,3,4,5,1,3,0,7,3,4,2,0,1,2,3,4,5,0,3,4,2,3,0,1,2,0,6,1,0,3,2,1,0,4,1,3,0,1,2,5,4,2,1,0,3,2,4,0,2,10,4,2,3,4,5,0,3,4,2,3,0,1,2,3,4,2,6,4,0,2,4,3,0,5,3,1,5,4,3,0,5,3,2,4,5,3,0,5,4,1,0,4,7,0,5,4,1,0,3,2,4,3,0,5,4,3,2,1,0,2,3,5,6,2,3,0,1,2,4,5,1,4,3,2,1,4,5,0,3,5,4,2,3,4,1,8,4,3,2,0,1,2,3,4,5,0,3,4,2,3,0,1,2,3,4,2,6,4,1,2,3,0,5,3,2,4,5,3,0,5,4,1,2,4,5,0,3,5,7,3,4,5,0,2,1,0,3,2,4,3,0,5,4,3,2,4,1,6,2,4,3,2,0,5,3,2,4,3,0,2,1,3,4,2,3,0,5,4,3,2,9,0,4,2,3,4,5,0,3,4,2,3,0,1,2,0,5,4,3,7,8,4,2,0,1,2,3,4,5,0,3,4,2,3,0,1,2,3,4,2,6,4,1,2,3,0,5,3,2,4,5,3,0,5,4,1,2,4,0,7,4,5,0,2,1,0,3,2,4,3,0,5,4,3,2,4,5,6,4,2,3,0,5,3,4,0,1,3,5,0,3,4,2,3,5,0,1,4,0,7,8,2,4,5,0,3,4,2,6,5,4,2,3,4,5,0,3,4,2,3,0,1,2,3,4,5,7,1,0,3,4,5,0,3,2,1,5,3,1,4,5,1,2,0,3,4,5,7,3,2,4,7,3,8,6,9,11,10,4,0,1,2,3,4,5,0,3,4,2,3,0,1,2,3,4,2,6,4,0,2,4,3,0,5,3,1,5,4,3,0,5,3,2,4,5,3,0,5,4,1,0,4,7,0,5,4,1,0,3,2,4,3,0,5,4,3,2,1,0,2,3,5,6,1,0,3,2,1,0,4,1,3,0,5,3,1,2,4,3,2,0,1,3,2,4,3,1,8,4,3,2,0,1,2,3,4,5,0,3,4,2,3,0,1,2,3,4,2,6,4,1,2,3,0,5,3,2,4,5,3,0,5,4,1,2,4,5,0,3,5,7,3,4,5,0,2,1,0,3,2,4,3,0,5,4,3,1,6,4,3,5,0,3,4,2,0,1,2,3,0,5,4,3,2,5,3,0,5,4,1,0,4,10,2,4,5,0,2,1,0,3,2,4,3,0,5,4,3,2,4,7,5,2,4,3,0,5,3,4,1,8,4,2,0,1,2,3,4,5,0,3,4,2,3,0,1,2,3,4,2,6,4,1,2,3,0,5,3,2,4,5,3,0,5,4,1,2,4,0,7,4,5,0,2,1,0,3,2,4,3,0,5,4,3,1,2,4,1,5,6,4,2,3,0,5,3,4,0,1,3,5,0,3,4,2,3,5,0,1,4,0,7,8,6,4,5,3,2,4,3,0,5,4,3,2,1,0,2,3,5,0,3,4,5,7,3,5,0,3,4,2,0,1,2,3,0,5,4,3,2,5,3,0,5,4,1,0,9,2,5,0,2,1,0,3,2,4,3,0,5,4,3,2,10,6,2,3,4,5,0,3,4,2,3,0,1,2,0,5,4,3,7,8,2,4,0,1,2,3,4,5,0,3,4,2,3,0,1,2,3,4,2,6,4,0,3,5,0,1,4,5,0,3,5,4,2,3,5,0,3,4,5,1,8,4,1,3,4,5,0,3,4,2,3,0,1,2,0,5,4,8,7,4,5,0,2,1,3,0,2,4,3,0,5,4,3,1,2,4,6,2,3,4,2,1,4,5,3,2,4,3,0,5,4,3,2,10,7,2,4,3,5,0,3,4,2,6,5,4,2,3,4,5,0,3,4,2,3,0,1,2,0,5,3,2,7,6,0,1,2,3,4,5,0,3,4,2,3,0,1,2,0,5,4,3,6,1,4,3,5,0,3,1,4,3,5,4,2,1,4,0,10,4,8,2,5,0,3,1,0,5,2,7,4,10,2,5,3,0,5,4,2,1,4,5,0,3,6,2,3,0,5,3,6,5,4,3,0,5,3,2,6,7,2,3,5,0,3,2,4,3,5,4,1,2,3,7,11,3,4,5,0,2,7,3,0,2,11,9,2,0,7,6,3,0,6,2,5,4,2,0,5,9,7,3,2,7,11,2,0,11,10,1,2,3,0,5,2,1,0,2,4,5,0,2,1,0,5,3,6,5,7,11,5,1,3,11,2,5,0,2,3,10,0,5,6,7,3,0,8,4,0,2,4,5,0,3,2,6,3,7,6,10,3,4,2,5,11}
\subsection{A Coil of length 2594 (Dimension: 13)}
\seqsplit{0,1,2,3,4,1,3,0,1,2,5,0,3,4,1,2,3,0,1,3,4,1,2,5,6,2,1,4,3,1,0,3,2,1,4,3,1,5,2,3,1,4,3,2,0,1,7,5,0,1,2,3,4,1,3,0,5,6,2,5,0,1,3,4,1,2,3,0,1,2,5,0,2,3,4,1,2,6,1,3,0,5,3,4,5,8,3,0,2,3,4,1,3,2,5,1,3,4,1,2,3,0,1,3,4,1,2,5,6,2,1,4,3,1,0,3,2,1,4,3,0,5,2,1,0,3,1,4,3,2,1,0,3,7,0,2,3,0,1,3,4,2,5,4,3,0,1,6,5,1,2,3,4,1,2,5,1,3,0,5,2,1,5,6,2,5,0,1,2,5,9,2,0,5,2,1,0,3,1,4,3,2,1,0,5,2,1,3,4,1,2,5,6,2,1,4,3,1,0,3,2,1,4,3,1,5,2,3,1,4,3,2,0,3,8,5,1,3,0,1,2,3,4,1,3,0,1,2,5,0,3,1,6,2,1,4,3,2,0,5,2,1,0,3,2,1,4,3,1,0,5,2,1,7,8,1,2,5,0,1,3,4,1,2,3,0,1,2,5,0,2,3,4,1,2,6,4,2,3,0,2,5,0,1,2,3,4,1,3,0,1,2,5,0,2,8,4,2,0,5,2,1,0,3,2,1,4,3,1,5,6,1,3,4,1,2,3,0,2,5,0,1,2,4,3,1,2,0,5,2,6,8,1,4,8,6,5,10,3,1,0,3,2,1,4,3,1,0,5,2,6,5,0,3,1,4,3,2,1,5,6,1,2,3,4,1,2,6,1,3,0,5,3,4,5,8,3,0,2,3,4,1,3,2,5,1,3,4,1,2,3,0,1,3,4,1,2,5,6,2,1,4,3,1,0,3,2,1,4,3,0,5,2,1,0,3,1,4,3,2,1,0,3,7,8,3,0,1,2,3,4,1,3,0,1,2,5,0,3,4,1,2,3,0,1,3,4,1,2,6,5,2,1,4,3,1,0,3,2,1,4,3,1,5,2,3,1,4,3,2,0,1,8,5,0,1,2,3,4,1,3,0,1,2,5,3,6,5,2,4,3,2,5,0,2,3,4,1,3,0,2,5,0,3,9,2,3,0,1,2,5,0,3,4,1,2,3,0,1,3,4,1,2,6,5,2,1,4,3,1,0,3,2,1,4,3,1,5,2,3,1,4,3,2,0,3,8,5,1,3,0,1,2,3,4,1,2,5,6,2,1,4,3,2,0,5,2,1,0,3,2,1,4,3,1,0,5,2,6,5,1,2,5,0,2,7,4,2,0,5,2,1,0,3,2,1,4,3,1,5,6,1,3,4,1,2,3,0,2,5,0,1,2,3,4,1,3,0,1,2,3,1,8,2,1,4,3,1,0,3,2,1,4,3,2,5,3,1,2,6,5,2,1,4,3,2,1,5,6,1,0,3,1,6,5,1,4,5,8,1,2,6,1,4,6,8,3,4,5,1,2,5,11,0,5,1,3,4,1,2,3,0,2,6,0,1,2,3,4,1,3,0,1,2,5,0,3,4,1,2,3,0,1,3,4,1,6,0,1,4,3,2,0,3,8,5,4,3,5,0,3,1,6,2,1,4,3,2,0,5,2,1,0,3,2,1,4,3,1,0,5,2,6,5,0,3,1,4,3,2,1,5,0,7,1,0,2,3,4,1,3,2,5,1,3,4,1,2,3,0,1,3,4,1,2,5,6,2,1,4,3,1,0,3,2,1,4,3,0,5,2,1,0,3,1,4,3,2,1,0,3,8,0,2,3,0,1,3,4,2,5,4,3,0,1,6,5,1,2,3,4,1,2,5,1,3,0,5,2,1,5,6,2,5,0,1,2,5,9,2,0,5,2,1,0,3,1,4,3,2,1,0,5,2,1,3,4,1,2,5,6,2,1,4,3,1,0,3,2,1,4,3,1,5,2,3,1,4,3,2,0,3,8,5,1,3,0,1,2,3,4,1,3,0,1,2,5,0,3,1,6,2,1,4,3,2,0,5,2,1,0,3,2,1,4,3,1,0,5,2,1,7,5,1,0,3,1,4,3,2,1,0,5,2,0,3,2,1,4,3,1,6,5,1,3,4,1,2,3,0,1,2,5,0,2,4,8,2,0,5,2,1,5,6,2,5,0,1,3,4,1,2,3,0,1,2,5,0,2,3,4,1,2,6,5,2,1,4,3,2,1,5,10,9,5,1,2,3,4,1,2,5,6,2,1,4,3,2,0,5,2,1,0,3,2,1,4,3,1,0,5,2,6,5,1,2,5,0,2,8,4,2,0,5,2,1,0,3,2,1,4,3,1,5,6,1,3,4,1,2,3,0,2,5,0,1,2,3,4,1,3,0,1,5,7,1,2,5,0,1,3,4,1,2,3,0,1,2,5,0,2,3,4,1,2,6,1,3,0,5,2,1,0,3,1,4,3,2,1,0,3,1,5,8,3,0,2,3,4,1,3,2,6,1,0,5,2,1,0,3,1,4,3,2,1,0,5,2,1,3,4,1,2,5,6,2,1,4,3,1,0,3,2,5,9,2,3,0,2,5,0,1,2,3,4,1,3,0,1,2,5,0,2,8,5,2,3,4,1,2,3,0,1,2,5,0,1,4,6,1,3,4,1,2,3,0,2,5,0,1,2,3,4,1,3,0,1,2,3,8,0,3,1,0,5,2,0,3,1,4,3,2,0,3,1,5,0,2,3,0,1,3,7,0,3,2,1,0,5,2,0,3,2,1,4,3,0,5,2,1,0,3,2,0,5,2,1,8,0,3,4,2,0,5,2,1,0,3,2,1,4,3,1,0,5,6,8,5,0,1,3,4,1,2,3,0,1,2,5,3,1,8,5,3,0,1,2,3,5,8,4,1,3,8,4,0,1,4,5,2,3,5,0,1,5,2,0,7,2,5,8,2,1,4,8,1,9,4,1,7,4,8,1,11,7,9,10,7,11,1,8,4,1,2,5,6,12,3,0,1,3,4,1,2,3,0,1,3,10,0,3,2,1,0,5,2,0,3,2,1,4,3,1,6,5,1,3,4,1,2,3,0,1,2,5,0,2,4,8,2,0,5,2,1,0,3,1,4,3,2,1,0,5,2,0,3,2,4,6,2,1,4,3,2,0,5,2,1,0,3,2,1,4,3,1,0,5,2,1,8,7,1,2,5,0,1,3,4,1,2,3,0,1,2,5,0,2,3,4,1,2,6,1,3,0,5,2,1,0,3,1,4,3,2,1,0,3,1,5,8,3,0,2,3,4,1,3,2,5,1,3,4,1,2,3,0,1,3,4,1,2,6,5,2,1,4,3,1,2,5,0,1,2,3,4,1,3,0,1,2,5,0,2,9,5,2,1,0,5,2,6,5,1,2,5,0,3,1,5,2,1,4,3,2,1,5,0,8,7,5,0,1,2,3,4,1,3,0,5,6,2,5,0,1,3,4,1,2,3,0,1,2,5,0,2,3,4,1,2,6,1,3,0,5,3,4,5,8,3,0,2,3,4,1,3,2,5,1,3,4,1,2,3,0,1,3,4,1,2,5,6,2,1,4,3,1,0,3,2,1,4,3,0,5,2,1,0,3,1,4,3,2,1,0,3,7,0,2,3,0,1,2,5,3,4,5,8,3,2,1,0,3,1,4,3,2,1,5,2,0,1,3,4,1,2,6,5,2,1,4,3,1,0,3,2,1,4,3,1,5,2,1,10,2,0,5,2,1,5,6,2,5,0,1,3,4,1,2,3,0,1,2,5,0,2,3,4,1,2,6,5,2,1,4,3,2,1,0,8,4,0,5,6,1,0,5,2,1,0,3,2,0,7,3,0,1,3,4,1,2,3,0,1,2,5,0,1,6,5,1,2,3,4,1,2,5,6,2,1,3,9,2,0,5,2,1,0,3,1,4,3,2,1,0,5,2,1,3,4,1,2,6,5,2,1,4,3,1,0,3,2,1,4,3,1,5,2,3,1,4,3,2,0,3,8,6,3,0,2,3,4,1,3,0,6,5,1,2,5,0,3,1,5,2,1,4,3,2,1,5,0,7,1,0,2,3,4,1,3,2,5,1,3,4,1,2,3,0,1,3,4,1,2,5,6,2,1,4,3,1,0,3,2,1,4,3,0,5,2,1,0,3,1,4,3,2,1,0,3,8,0,2,3,0,1,3,4,2,5,4,3,0,2,3,4,1,3,0,2,6,1,5,2,1,0,3,1,4,3,2,1,5,11,7,5,1,2,3,4,1,3,0,2,6,5,2,0,3,1,4,3,2,0,5,2,1,0,3,2,4,9,2,3,0,1,2,5,0,3,4,1,2,3,0,1,3,4,1,2,6,5,2,1,4,3,1,0,3,2,1,4,3,1,5,2,3,1,4,3,2,0,3,8,5,4,3,5,0,3,1,6,2,1,4,3,2,0,5,2,1,0,3,2,1,4,3,1,0,5,2,6,5,0,3,1,4,3,2,1,5,0,7,1,0,2,3,4,1,3,2,5,1,3,4,1,2,3,0,1,3,4,1,2,5,6,2,1,4,3,1,0,3,2,1,4,3,0,5,2,1,0,3,2,9,4,2,3,0,1,2,5,0,2,3,4,1,3,0,2,6,1,5,2,1,0,3,1,4,3,2,1,0,7,6,0,1,2,3,4,1,3,0,1,5,0,2,1,6,5,1,2,10,1,0,5,6,2,5,0,1,3,4,1,2,3,0,1,2,5,0,2,3,4,1,2,6,1,3,0,5,3,4,5,8,3,0,2,3,4,1,3,2,5,1,3,4,1,2,3,0,1,3,4,1,2,6,5,2,1,4,3,1,2,5,0,1,2,3,4,1,3,0,1,2,5,0,2,9,1,3,2,6,1,3,0,1,6,5,0,1,3,4,1,2,3,0,1,2,5,0,1,6,5,1,2,3,4,1,3,0,1,8,7,1,0,3,1,4,3,2,0,7,5,0,3,5,4,3,1,5,6,1,0,2,7,1,2,6,1,0,3,1,6,2,1,4,3,2,0,3,6,8,3,0,2,3,4,1,3,2,5,1,3,4,1,2,3,0,1,3,4,1,2,5,6,2,1,4,3,1,0,3,2,1,4,3,0,5,2,1,0,3,2,9,4,2,3,0,1,2,6,0,2,3,4,1,3,0,2,6,5,2,0,3,1,4,3,2,1,8,5,1,2,3,4,1,3,0,1,2,5,0,3,4,6,5,4,3,1,5,2,1,8,2,0,3,6,5,2,1,10,8,1,2,7,8,5,6,0,3,2,0,11,9,5,2,3,5,8,1,2,0,1,4,3,2,0,5,2,1,0,3,2,1,4,3,1,5,6,1,3,4,1,2,3,0,1,4,5,2,3,4,1,3,0,1,2,3,4,1,0,11,1,10,7,1,4,0,2,1,4,3,2,7,1,3,2,0,6,5,2,8,1,4,8,5,1,8,2,5,6,8,11,2,6,0,11,6,3,0,6,7,3,10,7,6,9,4,2,6,5,2,1,4,5,3,2,0,3,5,1,2,3,5,10,2,0,7,8,2,0,3,8,11,3,0,1,11,5,3,1,4,2,5,6,2,4,8,2,10,0,2,3,0,1,6,5,1,3,7,2,1,7,8,5,1,0,4,8,12,4,1,12,11,9,5,0,8,5,2,3,8,10,3,1,9,10,1,4,9,12,4,1,5,12,7,5,2,7,11,0,2,3,11,6,3,0,11,5,7,2,1,7,5,4,6,5,3,4,0,2,4,1,0,5,1,6,0,4,2,1,0,4,5,0,1,3,9,5,3,0,5,2,0,4,3,9,8,3,2,1,0,10}


\begin{thebibliography}{100}

\bibitem{AK91}
H. Abbott and M. Katchalski, On the construction of Snake In The Box codes, Utilitas Mathematica 40 (1991) 97--116.

\bibitem{AAC73}
L.E. Adelson, R. Alter and T.B, Curtz, Long snakes and a characterization of maximal snakes on the $d$-cube, Proc. 4th South-Eastern Conference on Combinatorics, Graph Theory and Computing, Congr. Numer. 8 (1973) 111--124.

\bibitem{Da65}
D.W. Davies, Longest ``separated'' paths and loops in an $n$ cube, IEEE Trans. Electronic Computers 14 (1965) 261.

\bibitem{Dr15}
T.E. Drapela, The Snake-in-the-Box Problem: a primer, MSc Thesis, Athens, Georgia, 2015.

\bibitem{EL00}
P.G. Emelyanova,  and A. Lukito, On the maximal length of a snake in hypercubes of small dimension,
Discrete Mathematics 218 (2000) 51--59.

\bibitem{GJ79}
M.R.Garey and D.S.Johnson, Computers and Intractability:A Guide to the Theory of \NP-completeness, Freeman, San Francisco, 1979.

\bibitem{HPB96}
A. P. Hiltgen, K. G. Paterson, and M. Brandestini. Single-track Gray codes, IEEE Transactions on Information Theory 42 (1996) 1555--1561.

\bibitem{HRSW15}
S. Hood, D. Recoskie, J. Sawada and D. Wong, Snakes, coils, and single-track circuit codes with spread $k$, Journal of Combinatorial Optimization 30 (2015) 42--62.

\bibitem{Ka58}
W.H. Kautz, Unit-distance error-checking codes, IRE Trans Electronic Computers EC-7 (1958) 179--180.

\bibitem{Ki12}
D. Kinny, A new approach to the Snake-In-The-Box-problem, Proc. ECAI 2012, 462--467.

\bibitem{Ko96}
K.J. Kochut, Snake-In-The-Box codes for dimension 7, Journal of Combinatorial Mathematics and Combinatorial Computing 20 (1996) 175--185.

\bibitem{MDWP15}
S.J. Meyerson, T.E. Drapela, W.E. Whiteside and W.D. Potter,
Finding longest paths in hypercubes: 11 new lower bounds for snakes, coils, and symmetrical coils, Proc. IEA/AIE 2015, 23--32.

\bibitem{MWDP14}
S.J. Meyerson, W.E. Whiteside, T.E. Drapela and W.D. Potter, Finding longest paths in hypercubes, snakes and coils, Proc. IEEE CIES 2014, 103--109.

\bibitem{OP15}
P.R.J. \"{O}sterg{\aa}rd and V.H. Pettersson, Exhaustive search for snake-in-the-box codes, Graphs and Combinatorics 31 (2015) 1019--1028.

\bibitem{OP14}
P.R.J. \"{O}sterg{\aa}rd and V.H. Pettersson, On the maximum length of coil-in-the-box codes in dimension~8, Discrete Applied Mathematics 179 (2014) 193--200.

\bibitem{PRMK94}
W. Potter, J. Robinson, J. Miller and K.J. Kochut, Using the genetic algorithm to find Snake- In-The-Box codes, Proc. IEA/AIE 1994, 421--426.

\bibitem{Si66}
R.C. Singleton, Generalized snake-in-the-box codes. IEEE Transactions on Electronic Computers (1966) 596--602.

\bibitem{Sn94}
H.S. Snevily, The snake-in-a-box problem: a new upper bound, Discrete Mathematics 133 (1994) 307--314.

\bibitem{Wo98}
  J. Wojciechowski, On constructing snakes in powers of complete graphs, Discrete Mathematics 181 (1998) 239--254

\bibitem{Wy12}
E. Wynn, Constructing circuit codes by permuting initial sequences, Manuscript, arXiv:1201.1647.

\bibitem{Ze97}
G. Z\'emor, An upper bound on the size of the snake-in-the-box, Combinatorica 17 (1997) 287--298.

 
\end{thebibliography}
\end{document}